\documentclass[journal,twoside]{IEEEtran}
\usepackage{cite}

\ifCLASSINFOpdf
 \usepackage[pdftex]{graphicx}
\else
\fi
\usepackage[cmex10]{amsmath}
\usepackage{array}

\usepackage{mdwmath}
\usepackage{mdwtab}
\usepackage[caption=false,font=footnotesize]{subfig}
\usepackage{stfloats}
\usepackage{url}

\usepackage[american,fulldiode]{circuitikz} 
\usepackage{graphicx} 
\usepackage{multirow}
\usepackage{amsfonts}
\usepackage{bbm}

\usepackage{multirow}
\usepackage[section]{placeins}
\usepackage{epstopdf}

\usepackage{slashbox}
\usepackage{amssymb}

\usepackage{lipsum,amsmath,multicol}

\usepackage{multirow}

\usepackage{float}
\usepackage{tablefootnote}
\usepackage[multiple]{footmisc}

\makeatletter
\def\ps@headings{%
\def\@oddhead{\mbox{}\scriptsize\rightmark \hfil \thepage}%
\def\@evenhead{\scriptsize\thepage \hfil \leftmark\mbox{}}%
\def\@oddfoot{}%
\def\@evenfoot{}}
\def\BState{\State\hskip-\ALG@thistlm}
\makeatother

\makeatletter
\newcommand{\removelatexerror}{\let\@latex@error\@gobble}
\makeatother

\usepackage{enumitem}

\allowdisplaybreaks

\usepackage{varwidth}
\newcolumntype{M}[1]{>{\begin{varwidth}[t]{#1}}l<{\end{varwidth}}}

\usepackage[linesnumbered,ruled,vlined]{algorithm2e}

\makeatletter
\newcommand{\algrule}[1][.2pt]{\par\vskip.5\baselineskip\hrule height #1\par\vskip.5\baselineskip}
\makeatother

\begin{document}
%
\title{A Method for Quickly Bounding the Optimal Objective Value of an OPF Problem using a Semidefinite Relaxation and a Local Solution}
%
%
%
\author{\IEEEauthorblockN{Alireza Barzegar,$^{\ast}$ \textit{Student Member,
      IEEE}}, \IEEEauthorblockN{Daniel K. Molzahn,$^{\dagger}$ \textit{Member,
      IEEE}, and \\ \IEEEauthorblockN{Rong Su,$^{\ast}$ \textit{Senior Member,
      IEEE}}}%
\thanks{\hspace*{-15pt}$\ast$: School of Electrical and Electronic Engineering, Nanayang Technological University, Singapore, alireza001@e.ntu.edu.sg, rsu@ntu.edu.sg.}%
\thanks{\hspace*{-15pt}$\dagger$: School of Electrical and Computer Engineering, Georgia Institute of Technology, molzahn@ece.gatech.edu.}}%

\maketitle

\begin{abstract}
Optimal power flow (OPF) is an important problem in the operation of electric power systems. Due to the OPF problem's non-convexity, there may exist multiple local optima. Certifiably obtaining the global optimum is important for certain applications of OPF problems. Many global optimization techniques compute an optimality gap that compares the achievable objective value corresponding to the feasible point from a local solution algorithm with the objective value bound from a convex relaxation technique. Rather than the traditional practice of completely separating the local solution and convex relaxation computations, this paper proposes a method that exploits information from a local solution to speed the computation of an objective value bound using a semidefinite programming (SDP) relaxation. The improvement in computational tractability comes with the trade-off of reduced tightness for the resulting objective value bound. Numerical experiments illustrate this trade-off, with the proposed method being faster but weaker than the SDP relaxation and slower but tighter than second-order cone programming (SOCP) and quadratic convex (QC) relaxations for many large test cases.
\end{abstract}



%
\IEEEpeerreviewmaketitle

\vspace*{-6pt}
\section{Introduction}
\label{l:introduction}
%
%
%
%

\IEEEPARstart{F}{irst} formulated by Carpentier~\cite{carpentier} in 1962, the optimal power flow (OPF) problem determines a steady-state operating point for an electric power system by optimizing an objective function, such as generation cost, while satisfying constraints that model the nonlinear power flow equations as well as limits on line-flows, voltage magnitudes, and generator outputs. The OPF problem is non-convex, NP-Hard~\cite{bienstock2015strong,pascalNPhard}, and may have multiple \emph{local optima} (i.e., feasible points that are superior to all nearby feasible points but possibly inferior to the \emph{global optimum})~\cite{bukhsh2013local}.

Research conducted for over fifty years has developed many mature local solution algorithms for OPF problems, including Newton-based techniques, sequential quadratic programming algorithms, interior point methods, etc.~\cite{opf_litreview1993IandII,ferc4,zimmerman2011matpower,frank2012optimal1,frank2012optimal2}. These algorithms often quickly find local optima for large OPF problems.

Certifiably finding global solutions to OPF problems (or, at least, strong objective value bounds) is important for many applications.
To compute global optimality certificates, recent research has focused on convex relaxations of OPF problems~\cite{molzahn_hiskens-fnt2017}. Convex relaxations bound the objective value of the original OPF problem, providing lower (upper) bounds for minimization (maximization) problems. Comparing the bound provided by a relaxation and the objective value obtained from a local solver yields an \emph{optimality gap}. A sufficiently small optimality gap certifies global optimality of the corresponding local solution.

Many OPF problems seek to minimize generation cost. Local optima in generation-cost-minimizing OPF problems may be significantly more expensive than the global solutions~\cite{bukhsh2013local}, and thus provably computing global optima is relevant to improving power system economics. Algorithms that use spatial branch-and-bound techniques to globally solve OPF problems are proposed in, e.g.,~\cite{gopalakrishnan2012,phan2012,harsha2018}. These algorithms iteratively apply local solvers and convex relaxation techniques in order to shrink the optimality gap until global optimality is certifiably achieved.

Global optimality guarantees are also important to OPF problems that are solved as subroutines in certain bi-level optimization algorithms, e.g., the robust OPF algorithm in~\cite{molzahn_roald-acopf_robust2018}. Examples of other applications where rigorous objective value bounds are particularly relevant include proving insolvability of power flow problems~\cite{molzahn2013sufficient}, bounding the errors associated with various power flow approximations~\cite{dvijotham2016error}, certifying security in the face of power injection variability for distribution systems with limited measurement and control capabilities~\cite{molzahn2019hicss}, and computing the feasible spaces of OPF problems for visualization purposes~\cite{molzahn2017computing}.



Various convex relaxation techniques provide a range of trade-offs between bound tightness and computational tractability. Simple relaxations, such as copper plate and network flow formulations~\cite{coffrin2016pscc}, quickly provide crude bounds. More sophisticated approaches, such as those based on second-order cone programming (SOCP)~\cite{jabr2006radial,low_tutorial1,sun2015,coffrin2016qc} and semidefinite programming (SDP)~\cite{lavaei2012zero,josz_molzahn-complex_hierarchy}, provide tighter bounds with associated trade-offs in computational speed. The literature details a variety of approaches for refining trade-offs in computational speed and tightness of the relaxations (e.g., adding valid cuts~\cite{sun2015} and exploiting problem structure~\cite{josz_molzahn-complex_hierarchy}). Also note that relaxation techniques are applicable to many types of systems, including mixed AC-DC networks~\cite{baradar2013second,bahrami2016semidefinite,venzke2019convex}. A comprehensive survey of relaxation techniques is provided in~\cite{molzahn_hiskens-fnt2017}. 

Building on the convex relaxation literature, this paper proposes a new method for balancing trade-offs between computational tractability and tightness in the calculation of objective value bounds for OPF problems. To speed computation of an objective value bound, our method takes as input an OPF solution obtained from a local solver. Our method therefore does not attempt to compute values for the decision variables; we emphasize that only the objective value bound is of interest in this paper. The proposed method is applicable as a subroutine in the bounding steps of global solution algorithms, such as those in~\cite{gopalakrishnan2012,phan2012,harsha2018}, as well as in other algorithms, such as those in~\cite{molzahn_roald-acopf_robust2018,molzahn2013sufficient,dvijotham2016error,molzahn2019hicss,molzahn2017computing}, that also require rigorous bounds on certain quantities. Reducing the time required to calculate objective value bounds thus has the potential to improve computational tractability for a variety of power system algorithms.

We specifically consider a method for obtaining an objective value bound using the SDP relaxation in~\cite{lavaei2012zero} and the dual variables from a local solution. The SDP relaxation in~\cite{lavaei2012zero} is derived from an OPF formulation where all non-convexity is contained in a rank condition. Following Shor's approach~\cite{shor1987}, an SDP relaxation is formed by neglecting this rank condition. A relaxation whose solution satisfies the rank condition is ``exact'', thus providing a tight bound on the objective value and the globally optimal decision variables. While the SDP relaxation in~\cite{lavaei2012zero} is exact for some OPF test cases, the rank condition is not always satisfied~\cite{lesieutre2011examining}. Nevertheless, the objective value from the SDP relaxation always bounds the objective of the original OPF problem.

Solving the SDP relaxation can be computationally challenging due to the presence of a positive semidefinite constraint on a large matrix. A matrix completion technique mitigates this computational challenge by decomposing the single positive semidefinite constraint into constraints on many submatrices~\cite{jabr2012exploiting,molzahn2013implementation}. However, the resulting formulation is still less tractable than many local solution algorithms. This motivates the development of techniques for exploiting information in a local solution to further speed computations of SDP relaxations.

Related work~\cite{molzahn2014sufficient} proposes a quickly checkable sufficient condition for global optimality of a candidate OPF solution from a local solver. Satisfaction of the Karush-Kuhn-Tucker (KKT) optimality conditions for the SDP relaxation certifies that a local solution is, in fact, globally optimal. Since the KKT conditions for the SDP relaxation can be quickly evaluated for a candidate OPF solution without explicitly solving an SDP, the condition in~\cite{molzahn2014sufficient} provides a tractable method for certifying global optimality.

While the approach in~\cite{molzahn2014sufficient} provides significant computational advantages for certain problems, the condition is sufficient but not necessary for global optimality and thus does not provide any information regarding global optimality when the condition is not satisfied. This paper proposes another method for leveraging information from a local solution to speed the computation of an objective value bound. Specifically, we seek a dual feasible point for the SDP relaxation. The objective value of any dual feasible point bounds the objective value of the primal SDP relaxation~\cite{boyd2009}, which itself bounds the objective value of the original OPF problem. Rather than attempt the computationally difficult task of computing the best achievable bounds, we instead leverage information from local solutions to more tractably obtain (possibly weaker) bounds. In particular, we fix certain dual variables in the SDP relaxation to their values from a local solution and solve the resulting simpler SDP problem.

Application to a variety of large test cases empirically demonstrates that the proposed method is faster than the SDP relaxation of~\cite{lavaei2012zero}, often without excessive degradation in the quality of the objective value bound. Conversely, empirical results indicate that the proposed method is slower than the SOCP relaxation in~\cite{jabr2006radial} and the QC relaxation in~\cite{coffrin2016qc} but can provide significantly tighter bounds. Thus, the proposed method can be viewed as a ``middle ground'' between the SDP relaxation and the SOCP and QC relaxations in terms of the trade-off between computational speed and tightness.

To summarize, the main contributions of this paper are the proposal and empirical analysis of a method for bounding the optimal objective value of an OPF problem by combining information from a locally optimal solution with SDP relaxation techniques. While the SDP relaxation techniques themselves are well known~\cite{lavaei2012zero,jabr2006radial,molzahn2013implementation}, the method proposed in this paper joins a very limited literature regarding the use of information from a local solver to rigorously obtain global information (i.e., objective value bounds) for OPF problems.

This paper is organized as follows. Section~\ref{l:opf} reviews formulations of the OPF problem and the SDP relaxation as well as the matrix completion decomposition. Section~\ref{proposed_Alg} proposes our method for quickly computing objective value bounds using the SDP relaxation in combination with a local solution, which is the main contribution of this paper. Section~\ref{Simulation_Results} numerically demonstrates this method using a variety of large test cases. Section~\ref{l:conclusion} concludes the paper.

\section{The OPF Problem and an SDP Relaxation}
\label{l:opf}
This section first formulates the OPF problem as well as its SDP relaxation~\cite{lavaei2012zero} and then discusses the matrix completion decomposition used to exploit network sparsity~\cite{jabr2012exploiting}.


\subsection{Optimal Power Flow Formulation}
Consider an $n$-bus system with the set of buses $\mathcal{N} = \left\lbrace 1, \ldots, n \right\rbrace$. Each bus~$k\in\mathcal{N}$ has an associated complex voltage phasor $V_{k}=V_{dk} + {j} V_{qk}$, power generation $P_{Gk} + {j} Q_{Gk}$, and specified power demand $P_{Dk} + {j} Q_{Dk}$, where $j = \sqrt{-1}$. The bus which sets the reference angle is denoted with the subscript ``ref''. The set of lines is denoted as $\mathcal{L}$, where $\left(l,m\right)\in\mathcal{L}$ indicates the line connecting buses~$l$ and~$m$. Each line $\left(l,m\right)\in\mathcal{L}$ is modeled as a $\Pi$ circuit with mutual admittance $y_{lm}$ and shunt susceptance $b_{sh,lm}$.\footnote{Extending our method to more general line models is straightforward. The numerical results in Section~\ref{Simulation_Results} use M{\sc atpower}'s line model~\cite{zimmerman2011matpower}.} The complex power flow and angle difference on each line $\left(l,m\right)\in\mathcal{L}$ are denoted by $S_{lm}=P_{lm}+jQ_{lm}$ and $\theta_{lm} = \theta_l - \theta_m$, respectively. The network admittance matrix is denoted as $\mathbf{Y}=\mathbf{G}+{j}\mathbf{B}$. Superscripts ``max'' and ``min'' indicate specified upper and lower limits. Buses without generators have the corresponding generation limits set to zero.



Denote $e_k$ as the $k^{th}$ standard basis vector in $\mathbb{R}^n$. Define the matrix $Y_k = e_{k}^{\vphantom{\intercal}} e_{k}^{\intercal}\mathbf{Y}$, where $\left(\,\cdot\,\right)^\intercal$ denotes the matrix transpose. Following the notation in~\cite{lavaei2012zero}, define the matrices $\mathbf{Y}_k$, $\mathbf{\overline{Y}}_k$, $\mathbf{M}_k$, and $\mathbf{N}_k$ used to construct expressions for the active and reactive power injections, the squared voltage magnitudes, and the angle reference, respectively, at bus~$k\in\mathcal{N}$:
\begin{subequations}
\begin{align}
\small
\label{eq:opf_Yk}
\mathbf{Y}_k & = \frac{1}{2}
\begin{bmatrix} \text{Re}( Y_{k} + Y_{k}^\intercal ) & \text{Im}( Y_{k}^\intercal - Y_{k} )  \\
\text{Im}( Y_{k} - Y_{k}^\intercal ) & \text{Re}( Y_{k} + Y_{k}^\intercal )
\end{bmatrix},\\
\label{eq:opf_Yk_}
\mathbf{\overline{Y}}_k & = -\frac{1}{2}
\begin{bmatrix} \text{Im}( Y_{k} + Y_{k}^\intercal ) & \text{Re}( Y_{k} - Y_{k}^\intercal )  \\
\text{Re}( Y_{k}^\intercal - Y_{k} ) & \text{Im}( Y_{k} + Y_{k}^\intercal )
\end{bmatrix}, \\
\label{eq:opf_Mk}
\mathbf{M}_k & =
\begin{bmatrix} e_{k}^{\vphantom{\intercal}} e_{k}^\intercal & 0  \\
\ 0 & e_{k}^{\vphantom{\intercal}} e_{k}^\intercal
\end{bmatrix}, \\
\label{eq:opf_Nk}
\mathbf{N}_k & =
\begin{bmatrix} 0 & 0  \\
\ 0 & e_{k}^{\vphantom{\intercal}} e_{k}^\intercal
\end{bmatrix}.
\end{align}
\end{subequations}
%


For each line~$\left(l,m\right)\in\mathcal{L}$, let $Y_{lm}=\left (\frac{b_{sh,lm}}{2}+y_{lm}\right )e_l^{\vphantom{\intercal}} e_l^\intercal - y_{lm} e_l^{\vphantom{\intercal}} e_m^\intercal$. To formulate the expressions for active and reactive power flows and phase angle differences associated with the line $\left(l,m\right)\in\mathcal{L}$, define the matrices
\begin{subequations}
\begin{align}
\small
\label{eq:opf_Ylm}
\mathbf{Y}_{lm} & =\frac{1}{2}
\begin{bmatrix} \text{Re}( Y_{lm} + Y_{lm}^\intercal ) & \text{Im}( Y_{lm} - Y_{lm}^\intercal )  \\
\text{Im}( Y_{lm}^\intercal - Y_{lm} ) & \text{Re}( Y_{lm} + Y_{lm}^\intercal )
\end{bmatrix},\\
\label{eq:opf_Ylm_}
\mathbf{\overline{Y}}_{lm} & =-\frac{1}{2}
\begin{bmatrix} \text{Im}( Y_{lm} + Y_{lm}^\intercal ) & \text{Re}( Y_{lm} - Y_{lm}^\intercal )  \\
\text{Re}( Y_{lm}^\intercal - Y_{lm} ) & \text{Im}( Y_{lm} + Y_{lm}^\intercal )
\end{bmatrix},\\
\label{eq:opf_Mlm}
\mathbf{M}_{lm} & =\frac{1}{2}
\begin{bmatrix} e_{l}^{\vphantom{\intercal}} e_{m}^\intercal + e_{m}^{\vphantom{\intercal}} e_{l}^\intercal & 0  \\
\ 0 & e_{l}^{\vphantom{\intercal}} e_{m}^\intercal + e_{m}^{\vphantom{\intercal}} e_{l}^\intercal
\end{bmatrix}, \\
\label{eq:opf_Mlm_}
\mathbf{\overline{M}}_{lm} & =\frac{1}{2}
\begin{bmatrix} 0 & -e_{l}^{\vphantom{\intercal}} e_{m}^\intercal + e_{m}^{\vphantom{\intercal}} e_{l}^\intercal  \\
\ e_{l}^{\vphantom{\intercal}} e_{m}^\intercal - e_{m}^{\vphantom{\intercal}} e_{l}^\intercal & 0
\end{bmatrix}.
\end{align}
\end{subequations}

By denoting the vector of voltage coordinates $x$ as
\begin{equation}
\small
\label{eq:opf_x vector}
x=\begin{bmatrix}
 V_{d1} &  V_{d2} &  \ldots &  V_{dn} &  V_{q1} &  V_{q2} &  \ldots &  V_{qn}
\end{bmatrix}
\end{equation}
and defining the rank-one matrix
\begin{equation}
\small
\label{eq:opf_W}
\mathbf{W}=xx^\intercal,
\end{equation}
all quantities of interest can be written in terms of $\mathbf{W}$. Let $ \mathrm{tr}\left( \cdot \right) $ denote the trace of a matrix. The active and reactive power injections at each bus~$k$ are $\mathrm{tr}\left( \mathbf{Y}_k \mathbf{W} \right)$ and $ \mathrm{tr}\left( \mathbf{\overline{Y}}_k \mathbf{W} \right)$, respectively. The squared voltage magnitude at bus~$k$ is $ \mathrm{tr}\left( \mathbf{M}_k \mathbf{W} \right)$. The active and reactive power flows into terminal~$l$ of line $\left(l,m\right)\in\mathcal{L}$ are $\mathrm{tr}\left( \mathbf{Y}_{lm} \mathbf{W} \right)$ and $\mathrm{tr}\left( \mathbf{\overline{Y}}_{lm} \mathbf{W} \right)$, and the phase angle difference is $\arctan \left( {\mathrm{tr}\left( \mathbf{\overline{M}}_{lm} \mathbf{W}\right)} / {\mathrm{tr}\left( \mathbf{M}_{lm} \mathbf{W}\right)}\right)$.

We consider a convex quadratic cost function of active power generation with specified quadratic, linear, and constant coefficients $c_{k2} \geq 0$, $c_{k1}$, and $c_{k0}$ for each generator $k\in\mathcal{N}$. Buses without generators have $c_{k2} = c_{k1} = c_{k0} = 0$.


Using these definitions, the OPF problem is

\begin{subequations}
\small
\vspace{-0.5em}
\label{eq:opf}
\begin{align}
\label{eq:opf_obj}
& \min \quad \sum_{k\in\mathcal{N}} \alpha_{k} \\
\nonumber
&\text{subject to}\quad (\, \forall k\in\mathcal{N},\; \left(l,m\right)\in\mathcal{L}\,) \\
\label{eq:opf_P}
&\quad  P_{Gk}^{min}- P_{Dk}\leq \mathrm{tr}\left( \mathbf{Y}_k \mathbf{W} \right)  \leq P_{Gk}^{max}- P_{Dk},\\
\label{eq:opf_Q}
&\quad  Q_{Gk}^{min}- Q_{Dk} \leq \mathrm{tr}\left( \mathbf{\overline{Y}}_k \mathbf{W} \right)  \leq Q_{Gk}^{max}- Q_{Dk},\\
\label{eq:opf_V2}
&\quad  \big(V_k^{min}\big)^2 \leq \mathrm{tr}\left( \mathbf{M}_k \mathbf{W} \right) \leq \big(V_k^{max}\big)^2,\\
\label{eq:opf_ref}
&\quad  \mathrm{tr}\left( \mathbf{N}_{ref} \mathbf{W} \right) = 0,\\
\label{eq:opf_S}
&\quad \begin{bmatrix}
 -\left( S_{lm}^{max} \right)^2 & \mathrm{tr}\left( \mathbf{Y}_{lm} \mathbf{W} \right) & \mathrm{tr}\left( \mathbf{\overline{Y}}_{lm} \mathbf{W} \right) \\
 \mathrm{tr}\left( \mathbf{Y}_{lm} \mathbf{W} \right) & -1 & 0 \\
 \mathrm{tr}\left( \mathbf{\overline{Y}}_{lm} \mathbf{W} \right) & 0 & -1
\end{bmatrix} \preceq 0, \\
\nonumber
&\quad \mathrm{tan} (\theta_{lm}^{min})\cdot\mathrm{tr}\left( \mathbf{M}_{lm} \mathbf{W} \right)  \leq \mathrm{tr}\left( \mathbf{\overline{M}}_{lm} \mathbf{W} \right)\\\label{eq:opf_PF_Eq_Phase_Angle}
&\quad \qquad\qquad\qquad\qquad\qquad\quad\;\; \ \leq  \mathrm{tan} (\theta_{lm}^{max})\cdot\mathrm{tr}\left( \mathbf{M}_{lm} \mathbf{W} \right), \\
\label{eq:opf_PF_Eq_P}
&\quad \begin{bmatrix}
c_{k1} \mathrm{tr}\left( \mathbf{Y}_k \mathbf{W} \right)  + a_{k} -\alpha_k &  \sqrt {c_{k2}} \mathrm{tr}\left( \mathbf{Y}_k \mathbf{W} \right) + b_{k}  \\
\sqrt {c_{k2}} \mathrm{tr}\left( \mathbf{Y}_k \mathbf{W} \right) + b_{k}  & -1
\end{bmatrix} \preceq 0, \\
\label{eq:opf_W_xxT}
&\quad  \mathbf{W}=xx^\intercal,
\end{align}
\end{subequations}
where $a_k = c_{k1}P_{Dk} + c_{k0}$, $b_k = \sqrt {c_{k2}}P_{Dk}$, and $\left(\,\cdot\,\right)\succeq 0$ indicates positive semidefinteness. Constraints~\eqref{eq:opf_P}--\eqref{eq:opf_V2} limit the power injections and squared voltage magnitudes, \eqref{eq:opf_ref} sets the reference angle, and~\eqref{eq:opf_S} and~\eqref{eq:opf_PF_Eq_Phase_Angle} limit the apparent power flows and phase angle differences. The objective~\eqref{eq:opf_obj} and constraint~\eqref{eq:opf_PF_Eq_P} represent the generation cost using a Schur complement formulation with auxiliary variables $\alpha_{k}$. All non-convexity is contained in the rank-one condition~\eqref{eq:opf_W_xxT}.

\subsection{Semidefinite Programming Relaxation of the OPF Problem}
\label{l:sdp_relax}

The approach in~\cite{lavaei2012zero} forms a Shor relaxation~\cite{shor1987} of~\eqref{eq:opf} by replacing the rank-one condition~\eqref{eq:opf_W_xxT} with a less stringent positive semidefinite matrix constraint,
\begin{equation}
\label{eq:sdpopf_W}
\mathbf{W}\succeq 0.
\end{equation}
This yields the primal formulation of the SDP relaxation:
%
\vspace{0.1em}
\begin{align}
\label{eq:sdp_primal}
& \min \quad \sum_{k\in\mathcal{N}} \alpha_{k} \quad \text{subject to\quad \eqref{eq:opf_P}--\eqref{eq:opf_PF_Eq_P},\; \eqref{eq:sdpopf_W}.}
\end{align}
%
If the solution to the SDP relaxation, $\mathbf{W}^\star$, satisfies the rank condition $\mathrm{rank}\left(\mathbf{W}^\star\right) = 1$, then the relaxation is \emph{exact} and thus provides both the globally optimal objective value and decision variables $V^\star = \sqrt{\lambda}\,\left(\eta_{1:n} + j\eta_{n+1:2n} \right)$, where $\lambda$ denotes the non-zero eigenvalue of $\mathbf{W}^\star$ with corresponding unit-length eigenvector $\eta$. While solutions for which $\mathrm{rank}\left(\mathbf{W}^\star\right) > 1$ do not provide globally optimal decision variables, the optimal objective value from the solution to~\eqref{eq:sdp_primal} still lower bounds the objective value of the original OPF problem~\eqref{eq:opf}.

Our method for speeding computations is based on the dual form of the SDP relaxation. To formulate the dual of~\eqref{eq:sdp_primal}, we define (for each bus~$k\in\mathcal{N}$ and each line $\left(l,m\right)\in\mathcal{L}$) the dual scalar variables ${\underline{\lambda}}_k$, ${\underline{\gamma}}_k$, ${\underline{\mu}}_k$, and ${\underline{\beta}}_{lm}$ associated with the lower bounds of \eqref{eq:opf_P}, \eqref{eq:opf_Q}, \eqref{eq:opf_V2} and \eqref{eq:opf_PF_Eq_Phase_Angle}; the dual scalar variables ${\overline{\lambda}}_k$, ${\overline{\gamma}}_k$, ${\overline{\mu}}_k$, and ${\overline{\beta}}_{lm}$ associated with the upper bounds of \eqref{eq:opf_P}, \eqref{eq:opf_Q}, \eqref{eq:opf_V2} and \eqref{eq:opf_PF_Eq_Phase_Angle}; the dual $3\times 3$ matrix variables $\mathbf{H}_{lm}$ associated with line-flow constraints \eqref{eq:opf_S}; and the dual $2\times 2$ matrix variables $\mathbf{R}_{k}$ associated with the Schur complement formulation of the quadratic cost function~\eqref{eq:opf_PF_Eq_P}.
%

For notational convenience, denote
%
\vspace{0.1em}
\label{eq:DUAL_Var}
\begin{align*}
& {\lambda}_k = {\overline{\lambda}}_k - {\underline{\lambda}}_k + c_{k1} + 2\sqrt{c_{k2}}R_k^{12},\quad {\gamma}_k = {\overline{\gamma}}_k - {\underline{\gamma}}_k, \\
& {\mu}_k = {\overline{\mu}}_k - {\underline{\mu}}_k, \quad \beta_{lm} = {\overline{\beta}}_{lm} - {\underline{\beta}}_{lm}.
\end{align*}
%
The dual formulation of the SDP relaxation~\eqref{eq:sdp_primal} is
\begin{subequations}
\label{eq:Dualsdpopf}
\begin{align}
& \max \quad \rho \\
\nonumber
&\text{subject to}\quad (\, \forall k\in\mathcal{N},\; \left(l,m\right)\in\mathcal{L}\,) \\
\label{eq:Dualsdpopf_A}
& \; \mathbf{A} \succeq 0, \\
\label{eq:Dualsdpopf_H}
& \;\mathbf{H}_{lm} \succeq 0,\\
\label{eq:Dualsdpopf_R}
& \;\mathbf{R}_{k} \succeq 0, \quad \mathbf{R}_{k}^{11}=1,\\
\label{eq:Dualsdpopf_Var}
&\; {\underline{\lambda}}_k \geq 0,\, {\overline{\lambda}}_k \geq 0,\, {\underline{\gamma}}_k \geq 0,\, {\overline{\gamma}}_k \geq 0,\, {\underline{\mu}}_k \geq 0,\, {\overline{\mu}}_k \geq 0, \\
&\; {\underline{\beta}}_{lm}\geq 0,\, {\overline{\beta}}_{lm}\geq 0,
\end{align}
\end{subequations}
where the scalar- and matrix-valued functions $\rho$ and $\mathbf{A}$ are
\small
\begin{align}
\nonumber 
& \rho = \sum_{k\in\mathcal{N}} \Bigg\{ \lambda_k P_{Dk} + \underline{\lambda}_k P_{Gk}^{min} - \overline{\lambda}_k P_{Gk}^{max}+\gamma_k Q_{Dk} + \underline{\gamma}_k Q_{Gk}^{min} \\ \nonumber
& \quad - \overline{\gamma}_k Q_{Gk}^{max}+ \underline{\mu}_k \big(V_k^{min}\big)^2 - \overline{\mu}_k \big(V_k^{max}\big)^2 \Bigg\} + \sum_{k\in\mathcal{N}} \Big( c_{k0} - \mathbf{R}_{k}^{22} \Big) \\\label{eq:DualCost_rho}
& \quad - \sum_{(l,m)\in\mathcal{L}} \Bigg\{ \big( S_{lm}^{max} \big)^2 \mathbf{H}_{lm}^{11} + \mathbf{H}_{lm}^{22} + \mathbf{H}_{lm}^{33} \Bigg\},\\
\nonumber
& \mathbf{A} = \sum_{k\in\mathcal{N}} \Bigg\{ \lambda_k \mathbf{Y}_{k} + \gamma_k \overline{\mathbf{Y}}_k + {\mu}_k \mathbf{M}_k \Bigg\}\\\nonumber
& \quad  + \sum_{(l,m)\in\mathcal{L}} \Bigg\{ 2\mathbf{H}_{lm}^{12} \mathbf{Y}_{lm} + 2\mathbf{H}_{lm}^{13} \overline{\mathbf{Y}}_{lm}^{13}  - \overline{\mathbf{M}}_{lm} \beta_{lm} \\\label{eq:DualMat_A}
& \qquad \qquad + \mathrm{tan}\Big(\theta_{lm}^{max}\Big)\mathbf{M}_{lm}{\overline{\beta}}_{lm} - \mathrm{tan}\Big(\theta_{lm}^{min}\Big)\mathbf{M}_{lm}{\underline{\beta}}_{lm} \Bigg\}.
\end{align}
\normalsize
Note that $\mathbf{R}_{k}^{cd}$ and $\mathbf{H}_{lm}^{cd}$ represent the $ (c,d)$ elements of matrices $\mathbf{R}_{k}$ and $\mathbf{H}_{lm}$, respectively.

Any feasible point for the dual SDP relaxation~\eqref{eq:Dualsdpopf} provides a lower bound for the objective value of the primal SDP relaxation~\eqref{eq:sdp_primal}~\cite{boyd2009}, and thus a lower bound for the original OPF problem~\eqref{eq:opf}. The solution to~\eqref{eq:Dualsdpopf} provides the best achievable lower bound among all dual feasible points.


\subsection{Matrix Completion Decomposition}
\label{l:matrix_completion}

The positive semidefinite constraint on a $2n\times 2n$ matrix ($\mathbf{W}\succeq 0$ in the primal form, $\mathbf{A}\succeq 0$ in the dual form) challenges the computational tractability of the SDP relaxation. To improve computational speed, a matrix completion theorem~\cite{gron1984} enables the decomposition of this single positive semidefinite constraint into constraints on many smaller submatrices~\cite{fukuda2001exploiting,jabr2012exploiting}. Our method leverages this decomposition.

This section summarizes the matrix completion decomposition, beginning with several definitions from graph theory. A ``clique'' is a set of nodes which are adjacent to all other nodes in the clique.\footnote{Radial sections of the network have associated maximal cliques with two buses, while cyclic portions of the network have associated maximal cliques with three or more buses.} A clique that is not contained in a larger clique is a ``maximal clique''. An edge that connects two non-adjacent vertices in a cycle is called a ``chord''. A graph is ``chordal'' if every cycle with the length of four or more nodes has a chord. We use $\mathcal{S}$ to denote the set of all maximal cliques and $m$ to indicate the number of maximal cliques.

\begin{figure}[t]
 \centering
 \includegraphics[width=6cm, height=2.5cm]{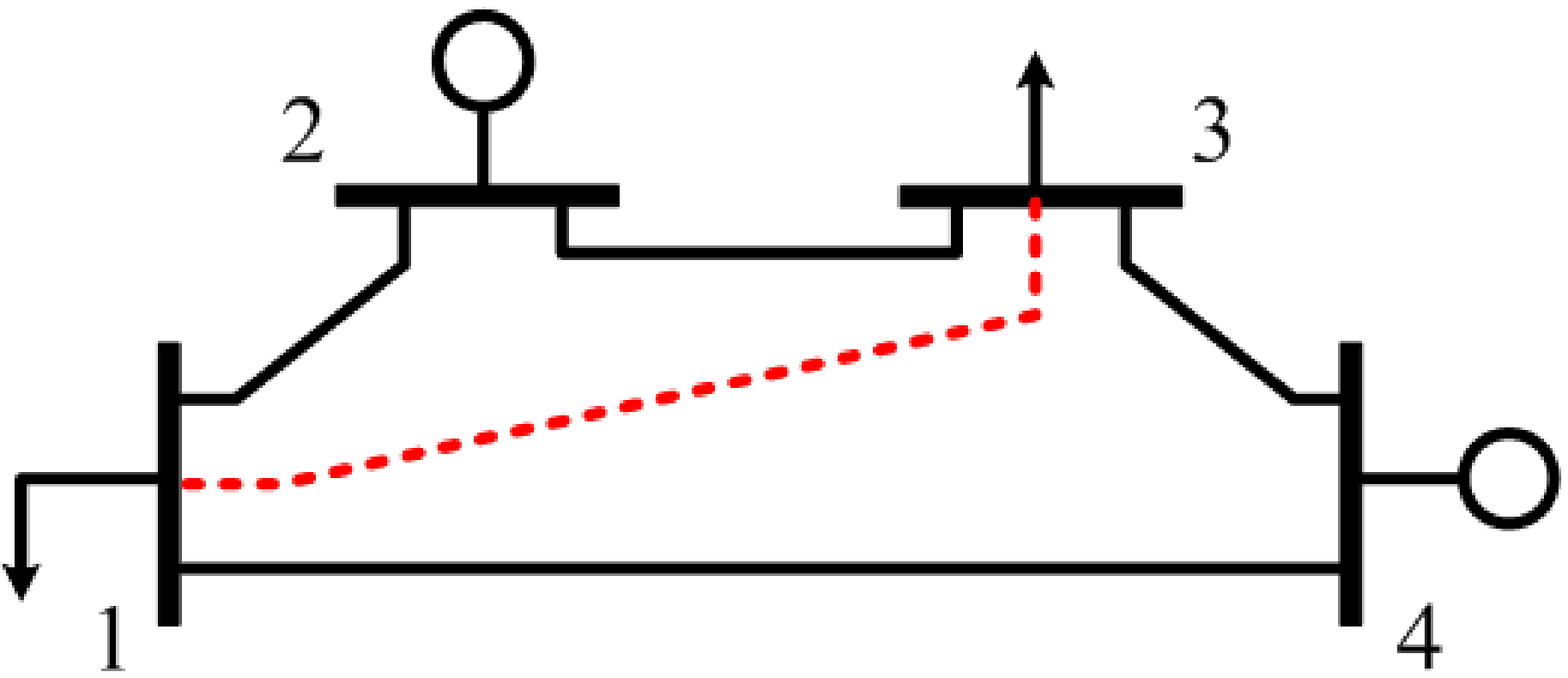}
 \caption{Example of the one-line diagram of a four-bus power network (black lines) and a chordal extension (black lines plus the dashed red line).}
 \label{One-line_diagram_of_a_4-bus_power_network}
 \vspace*{-1em}
 \end{figure}

As a simple example, consider the one-line diagram of a four-bus system shown in Fig.~\ref{One-line_diagram_of_a_4-bus_power_network}. The network graph is not chordal since the cycle $\{1, 2, 3, 4\}$ does not have a chord. The graph formed by adding the dashed red line is chordal with two maximal cliques $\mathcal{C}_1=\{1, 2, 3\}$ and $\mathcal{C}_2=\{1, 3, 4\}$. Thus, the graph including the dashed red line is a ``chordal extension''. See~\cite{jabr2012exploiting} and~\cite{valiente2013algorithms} for further details on these definitions.

The matrix completion theorem in~\cite{gron1984} considers a partially specified symmetric matrix with an associated undirected chordal graph. For our purposes, this graph is (a chordal extension of) the power system network and the unspecified matrix entries correspond to non-adjacent pairs of buses. The theorem states that the partially specified matrix can be ``completed'' to be positive semidefinite (i.e., there exist values for the unspecified entries such that the matrix is positive semidefinite) if and only if all the submatrices corresponding to the graph's maximal cliques are positive semidefinite.



For non-chordal power system networks, the matrix completion decomposition is applied by forming a chordal extension of the network graph. An appropriate chordal extension is constructed using the sparsity pattern of a Cholesky factorization of the network's adjacency matrix~\cite{fukuda2001exploiting,jabr2012exploiting}. A minimum degree ordering reduces the number of additional edges required in the chordal extension~\cite{amestoy2004algorithm}. Further computational improvements are possible via heuristics for constructing alternative chordal extensions~\cite{molzahn2013implementation}.

Using the matrix completion theorem, the positive semidefinite constraint is equivalently replaced by constraints on smaller submatrices associated with the maximal cliques of the chordal extension of the network graph. For the system in Fig.~\ref{One-line_diagram_of_a_4-bus_power_network}, this implies that the constraint $\mathbf{W} \succeq 0$ in~\eqref{eq:sdpopf_W} can be replaced by the constraints on the two submatrices whose rows and columns are associated with buses in the maximal cliques $\mathcal{C}_1$ and $\mathcal{C}_2$. Due to the possibility of non-empty intersections among maximal cliques, note that different submatrices may contain elements which refer to  the same element in the original $2n\times2n$ matrix $\mathbf{W}$. ``Linking constraints'' ensure consistency among variables contained in the submatrices associated with multiple cliques~\cite{fukuda2001exploiting,jabr2012exploiting}.


As an illustrative example, again consider the chordal extension of the four-bus system in Fig.~\ref{One-line_diagram_of_a_4-bus_power_network}. The matrix completion theorem enables the decomposition of the $8\times 8$ positive semidefinite matrix constraint~\eqref{eq:sdpopf_W} into two constraints on the $6\times 6$ submatrices, denoted $\mathbf{W}^{(1)}$ and $\mathbf{W}^{(2)}$, that are associated with the cliques $\mathcal{C}_1$ and $\mathcal{C}_2$. For notational simplicity, we next show the upper-left $3\times 3$ diagonal blocks and corresponding linking constraints for these submatrices. The remaining entries are treated analogously.

\vspace*{-1em}
\begin{subequations}
\small
\begin{equation*}
\label{eq:submatrices_W_}
\begin{array}{ll}
\mathbf{W}^{(1)} =  & \mathbf{W}^{(2)} = \\
\begin{bmatrix}
w^{(1)}_{11} & w^{(1)}_{12} & w^{(1)}_{13}  & \!\!\!\cdots \\
\ w^{(1)}_{12} & w^{(1)}_{22} & w^{(1)}_{23} & \!\!\!\cdots \\
\ w^{(1)}_{13} & w^{(1)}_{23} & w^{(1)}_{33} & \!\!\!\cdots \\[-0.25em]
\vdots & \vdots & \vdots & \!\!\!\ddots
\end{bmatrix} \succeq 0,
& 
\begin{bmatrix}
w^{(2)}_{11} & w^{(2)}_{12} & w^{(2)}_{13} & \!\!\!\cdots \\
\ w^{(2)}_{12} & w^{(2)}_{22} & w^{(2)}_{23} & \!\!\!\cdots \\
\ w^{(2)}_{13} & w^{(2)}_{23} & w^{(2)}_{33} & \!\!\!\cdots \\[-0.25em]
\vdots & \vdots & \vdots & \!\!\!\ddots
\end{bmatrix}  \succeq 0
\end{array}
\end{equation*}
\end{subequations}
Since buses $1$ and $3$ are the first and third elements in $\mathcal{C}_1$ and the first and second elements in $\mathcal{C}_2$, the linking constraints are
%
\begin{equation*}
\vspace{0.1em}
\label{eq:example_link_const}
\begin{split}
w^{(1)}_{11}=w^{(2)}_{11},\quad w^{(1)}_{13}=w^{(2)}_{12},\quad
& w^{(1)}_{33}=w^{(2)}_{22},\; \ldots
\end{split}
\end{equation*}

Applying this decomposition to the matrix $\mathbf{W}$ in~\eqref{eq:sdpopf_W} significantly improves the computational tractability of the SDP relaxation~\cite{jabr2012exploiting}. Taking the dual of the resulting SDP yields a chordal-sparsity-exploiting version of~\eqref{eq:Dualsdpopf}, which is the formulation used by our proposed method. We hereafter denote the submatrices in the sparsity-exploiting version of~\eqref{eq:Dualsdpopf} as $\mathbf{A}_i$, $i=1,\ldots,m$, each corresponding to a maximal clique contained in the set $\mathcal{S} = \left\lbrace \mathcal{C}_1,\ldots,\mathcal{C}_m\right\rbrace$. Note that the sparsity-exploiting version of~\eqref{eq:Dualsdpopf} includes dual variables in $\mathbf{A}_i$ that correspond to the linking constraints in the primal version. For further details, see~\cite{molzahn_hiskens-fnt2017,jabr2012exploiting,molzahn2013implementation}.

\section{Proposed Method}
\label{proposed_Alg}
Despite the computational improvements from the matrix completion decomposition, sparsity-exploiting SDP relaxations are slower than many local solution algorithms for typical OPF problems. This section describes our proposed method for leveraging the speed of local solution algorithms to more quickly compute an objective value bound using the SDP relaxation. We present a high-level summary, detail each step, and discuss several features of the proposed method.

\subsection{Description of the Proposed Method}
\label{Algo_descrip}
Our method is based on the fact that the objective value associated with any dual feasible point for the SDP relaxation~\eqref{eq:Dualsdpopf} bounds the objective of the primal SDP relaxation~\eqref{eq:sdp_primal} and therefore also bounds the objective of the original OPF problem~\eqref{eq:opf}. Rather than perform the computationally expensive calculation of the best achievable bound via application of an SDP solver to~\eqref{eq:Dualsdpopf}, we instead use certain information from a local solution to form a simplified SDP problem. Applying an SDP solver to this simplified problem can quickly provide a dual feasible point whose corresponding objective value is a bound on the objective of the original OPF problem~\eqref{eq:opf}. The trade-off inherent to this computational improvement is potential suboptimality in the obtained bound.



To form this simplified problem, we fix certain dual variables in the SDP relaxation~\eqref{eq:Dualsdpopf} to their corresponding values from a specified local solution, leaving only a subset of the dual variables to vary. The fixed dual variables are identified using a heuristic that computes ``problematic'' submatrices of the matrix $\mathbf{A}$ in~\eqref{eq:Dualsdpopf} evaluated using the dual variable values from the local solution. Consider the submatrices $\mathbf{A}_i$, $i=1,\ldots,m$, resulting from the matrix completion decomposition. A submatrix is deemed ``problematic'' if it has any eigenvalues that are either negative or small positive numbers. Dual variables associated with ``problematic'' submatrices are allowed to vary, while the other dual variables in~\eqref{eq:Dualsdpopf} are fixed to their values from the local solution.

Algorithm~\ref{alg:quick} summarizes our method. Each step is described in further detail below. 

\begin{figure}[t]
 \removelatexerror
  \begin{algorithm}[H]
\label{alg:quick}
\caption{Quickly Bounding the Objective Value}
\SetAlgoLined
\KwIn{Dual variables from a local OPF solution.}
\KwOut{A bound for the OPF's optimal objective value.}\vspace*{-3pt}\algrule[0.5pt]\vspace*{-3pt}
  Obtain the matrix $\hat{\mathbf{A}}$ matrix by substituting the dual variable values from the local solution into~\eqref{eq:DualMat_A}.\\
  Compute the matrix completion decomposition, i.e., the submatrices $\mathbf{A}_i$ for all $i\in\mathcal{S}$ using the approach in~\cite{jabr2012exploiting,molzahn2013implementation}. Substitute in the values of the dual variables from the local solution to form $\hat{\mathbf{A}}_i$.\\
  Identify the set of problematic submatrices, $\mathcal{S}_p$.\\
  Create the sets of all buses, $\mathcal{N}_{p}$, and lines, $\mathcal{L}_{p}$, included in $\mathcal{S}_{p}$.\\
  Formulate and solve the simplified SDP problem~\eqref{eq:NEW_Dual_sdpopf}. Return the objective value as a bound on the original OPF problem~\eqref{eq:opf}.\\
 \end{algorithm}
 \vspace*{-1em}
\end{figure}

\begin{enumerate}[wide]
 \item[\textbf{Input}:] Our method takes as input the values of the dual variables from a local solution to the OPF problem~\eqref{eq:opf}. These values can be obtained using any primal/dual local OPF solver.  
 
Note that solutions to equivalent but differently formulated OPF problems can have different values for the dual variables. For instance, the OPF formulation in M{\sc atpower}~\cite{zimmerman2011matpower} limits the voltage magnitudes rather than the squared voltage magnitudes as in~\eqref{eq:opf_V2}. If the local solution is computed for a different OPF formulation than~\eqref{eq:opf}, use conversions such as those in~\cite{molzahn2014sufficient} to obtain dual variables that are consistent with~\eqref{eq:opf}.
 
 
 \item[\textbf{1)}] \textbf{Evaluate the dual matrix}: Substitute the values of the dual variables into the $\mathbf{A}$ matrix defined in~\eqref{eq:DualMat_A}. Denote the resulting matrix as $\hat{\mathbf{A}}$.
 
 \item[\textbf{2)}] \textbf{Construct the matrix completion decomposition}: As summarized in Section~\ref{l:matrix_completion}, use the approach in~\cite{jabr2012exploiting,molzahn2013implementation} to compute the maximal cliques of a chordal extension of the network graph and form the corresponding submatrices of $\hat{\mathbf{A}}$, i.e., $\hat{\mathbf{A}}_1,\ldots,\hat{\mathbf{A}}_m$.\footnote{ While enumerating the maximal cliques of a general graph is NP-hard, the maximal cliques of a \emph{chordal} graph can be computed in linear time~\cite{tarjan1984simple}. Thus, Step~3 of our method is computationally tractable.} When constructing $\hat{\mathbf{A}}_i$, the dual variables corresponding to the linking constraints are set to zero.
 
 \item[\textbf{3)}] \textbf{Identify problematic submatrices}: Use the following heuristic to identify a set of ``problematic'' cliques, which is denoted as $\mathcal{S}_{p}$. Let $|\,\cdot\,|$ indicate the cardinality of a set. Let $\mathrm{ceil}(\,\cdot\,)$ denote the ceiling function which returns the smallest integer greater than or equal to the argument. Let $\mathrm{eig}(\,\cdot\,)$ return the eigenvalues of a matrix. Let $\mathcal{S}^{({\nsucceq 0})}$ define the set of cliques associated with the submatrices that are not positive semidefinite.
 
 Define a scalar parameter $\sigma$ such that $0 \leq \sigma \leq 1$. The parameter $\sigma$ dictates the minimum percentage of submatrices that are identified as ``problematic''. The set of problematic cliques, $\mathcal{S}_p$, consists of the cliques associated with the first $\mathrm{ceil}(\sigma\cdot|\mathcal{S}|)$ submatrices $\hat{\mathbf{A}}_i$ listed in order of their smallest eigenvalues, $\mathrm{min}(\mathrm{eig}(\hat{\mathbf{A}}_i))$. In order to ensure that the cliques associated with all non-positive-semidefinite submatrices are identified as ``problematic'', choose $\sigma$ such that  $\sigma \geq \left|\mathcal{S}^{({\nsucceq 0})}\right| / \left|\mathcal{S}\right|$.
 
 \item[\textbf{4)}] \textbf{Identify the fixed and non-fixed dual variables}: Identify the set $\mathcal{N}_p$ that contains all buses which are in problematic cliques $\mathcal{S}_p$. Similarly, let $\mathcal{L}_p$ denote the set of all lines for which either terminal bus is in a problematic clique $\mathcal{S}_p$. The dual variables associated with buses in $\mathcal{N}_p$ and lines in $\mathcal{L}_p$ are decision variables in the simplified SDP problem. Dual variables associated with the remaining buses and lines, i.e., $\mathcal{N}\setminus\mathcal{N}_p$ and $\mathcal{L}\setminus\mathcal{L}_p$, are fixed to their corresponding values from the local solution.
 
 \item[\textbf{5)}] \textbf{Formulate and solve the simplified SDP relaxation}: Fixing the dual variables for buses in $\mathcal{N}\setminus\mathcal{N}_p$ and lines in $\mathcal{L}\setminus\mathcal{L}_p$ yields the following simplified version of~\eqref{eq:Dualsdpopf}:
 \begin{subequations}
\label{eq:NEW_Dual_sdpopf}
\begin{align}
& \max \quad \rho \\
\nonumber
&\text{subject to} \quad \left(~\forall k\in\mathcal{N}_{p},~ \forall (l,m)\in\mathcal{L}_{p}\right)\hspace*{-90pt}\\
\label{eq:NEW_Dual_sdpopf_A}
& \; \tilde{\mathbf{A}}_i \succeq 0, & \forall i=1,\ldots,m\\
\label{eq:NEW_Dual_sdpopf_H}
& \; \mathbf{H}_{lm} \succeq 0, \\
\label{eq:NEW_Dual_sdpopf_R}
& \; \mathbf{R}_{k} \succeq 0, \quad \mathbf{R}_{k}^{11}=1, \\
\label{eq:NEW_Dualsdpopf_Var}
& \; {\underline{\lambda}}_k \geq 0,\, {\overline{\lambda}}_k \geq 0,\, {\underline{\gamma}}_k \geq 0,\, {\overline{\gamma}}_k \geq 0,\, {\underline{\mu}}_k \geq 0,\, {\overline{\mu}}_k \geq 0,\hspace*{-90pt}  \\
&\; {\underline{\beta}}_{lm}\geq 0,\, {\overline{\beta}}_{lm}\geq 0,
\end{align}
\end{subequations}
%
%
where the submatrices in the sparsity-exploiting dual formulation are denoted as $\tilde{\mathbf{A}}_i$. These submatrices include constant terms which correspond to the values of the dual variables associated with the buses in~$\mathcal{N}\setminus\mathcal{N}_p$ and lines in~$\mathcal{L}\setminus\mathcal{L}_p$, i.e., the fixed dual variable values from the local solution. Furthermore, all matrices $\tilde{\mathbf{A}}_i$, $i=1,\ldots,m$, are functions of the dual variables corresponding to the linking constraints, as discussed in Section~\ref{l:matrix_completion} and~\cite{ molzahn_hiskens-fnt2017,jabr2012exploiting,molzahn2013implementation}. The dual variables corresponding to the linking constraints are decision variables in~\eqref{eq:NEW_Dual_sdpopf} whose values are not provided by the local solution. Thus,~\eqref{eq:NEW_Dual_sdpopf_A} considers the submatrices $\tilde{\mathbf{A}}_i$ corresponding to \emph{all} maximal cliques, $\mathcal{S}$, rather than just those corresponding to the problematic cliques $\mathcal{S}_p$.

A solution to~\eqref{eq:NEW_Dual_sdpopf} is a dual feasible point for the (sparsity-exploiting) SDP relaxation~\eqref{eq:Dualsdpopf}. Thus, the objective value for the solution to~\eqref{eq:NEW_Dual_sdpopf} is a lower bound on the objective value of the original OPF problem~\eqref{eq:opf}.

\end{enumerate}
 
\subsection{Discussion}
\label{III:discustion}
With many variables fixed to their corresponding values from the local solution, the SDP problem~\eqref{eq:NEW_Dual_sdpopf} is easier to solve than the original SDP relaxation~\eqref{eq:Dualsdpopf}, which results in a computational speed advantage for the proposed method. This advantage comes at the cost of reduced tightness since the dual problem has fewer degrees of freedom with which to improve the objective $\rho$. The value of $\sigma$ controls this trade-off between the solution speed and bound tightness. The tightest lower bound occurs for $\sigma=100\%$, which yields the (sparsity-exploiting) complete dual SDP relaxation~\eqref{eq:Dualsdpopf}, i.e., all dual variables are allowed to vary. Conversely, $\sigma = 0\%$ fixes the largest number of dual variables in~\eqref{eq:NEW_Dual_sdpopf} and thus yields the fastest computational speed at the cost of a weaker bound. 

The proposed method performs well with small values of $\sigma$ when the local solution is close to the SDP solution. To numerically evaluate this claim, we compute the two-norm $\left|\left|\,\cdot\,\right|\right|_2$ of the ratio between the dual variable values for the local solution and the SDP solution:
\begin{equation}
\small
\label{eq:norm_Dist}
 \text{Dual Correspondence Ratio (\%)} = \left|\left|\frac{X_{local}^{Dual}-X_{SDP}^{Dual}}{X_{SDP}^{Dual}}\right|\right|_2 \times 100,
\end{equation}
where $X$ denotes the vector of the dual variables defined in where $X$ denotes the vector of the dual variables defined in Section~\ref{l:sdp_relax}. When the local solution is far from the SDP solution, fixing variables to the values from the local solution overly constrains~\eqref{eq:NEW_Dual_sdpopf}, resulting in an inferior bound and potentially even infeasibility. Such cases can be addressed by increasing the value of $\sigma$.

The dual correspondence ratio in~\eqref{eq:norm_Dist} provides an a-posterori explanation of our method's performance. We also consider a heuristic for estimating the quality of our method's bound without prior knowledge of the dual SDP solution. Specifically, we empirically find that our method performs well when more than $95$\% of the submatrices $\hat{\mathbf{A}}_i$ obtained using the local solution's dual variables are positive semidefinite. We suggest using this heuristic as a prescreening step to predict whether our method can be expected to perform well or, conversely, if alternative relaxation techniques should be applied.

\section{Numerical Experiments}
\label{Simulation_Results}
A key advantage of the proposed method is its ability to control the trade-off between computational speed and tightness of the objective value bound. This trade-off is determined by how well the dual variables from the local solution match the dual variables for the SDP relaxation's solution. The closeness of these solutions is problem dependent. Thus, evaluating the effectiveness of the proposed method requires numerical experimentation. After providing our implementation details and evaluation methodology, this section applies the proposed method to a variety of large test cases. The results indicate that the proposed method is a valuable ``middle ground'' in that it is often tighter but slower than the SOCP and QC relaxations while weaker but faster than the SDP relaxation.

\subsection{Implementation Details and Evaluation Methodology}
Algorithm~\ref{alg:quick} is implemented using MATLAB~R2016a and YALMIP~\cite{lofberg2004yalmip} and solved using MOSEK~8 on a computer with a 64-bit Intel Quad Core 3.40~Ghz CPU with 8~GB of RAM. The local solutions are obtained using the interior point method in M{\sc atpower}~\cite{zimmerman2011matpower}. The numerical results include comparisons to the SOCP relaxation~\cite{low_tutorial1}, the QC relaxation~\cite{coffrin2016qc} augmented with the ``Lifted Nonlinear Cuts'' from~\cite{coffrin2016strengthen_tps}, and the sparsity-exploiting SDP relaxation~\cite{jabr2012exploiting,molzahn2013implementation}. For our experiments, we consider the large test cases ($\geq 118$ buses) in the NESTA 0.6.0 archive~\cite{nesta}. The remainder of this section details results for a representative set of these test cases.

\begin{figure*}[!t]
\centering
\subfloat[case2746wop\_mp]{\includegraphics[totalheight=0.17\textheight]{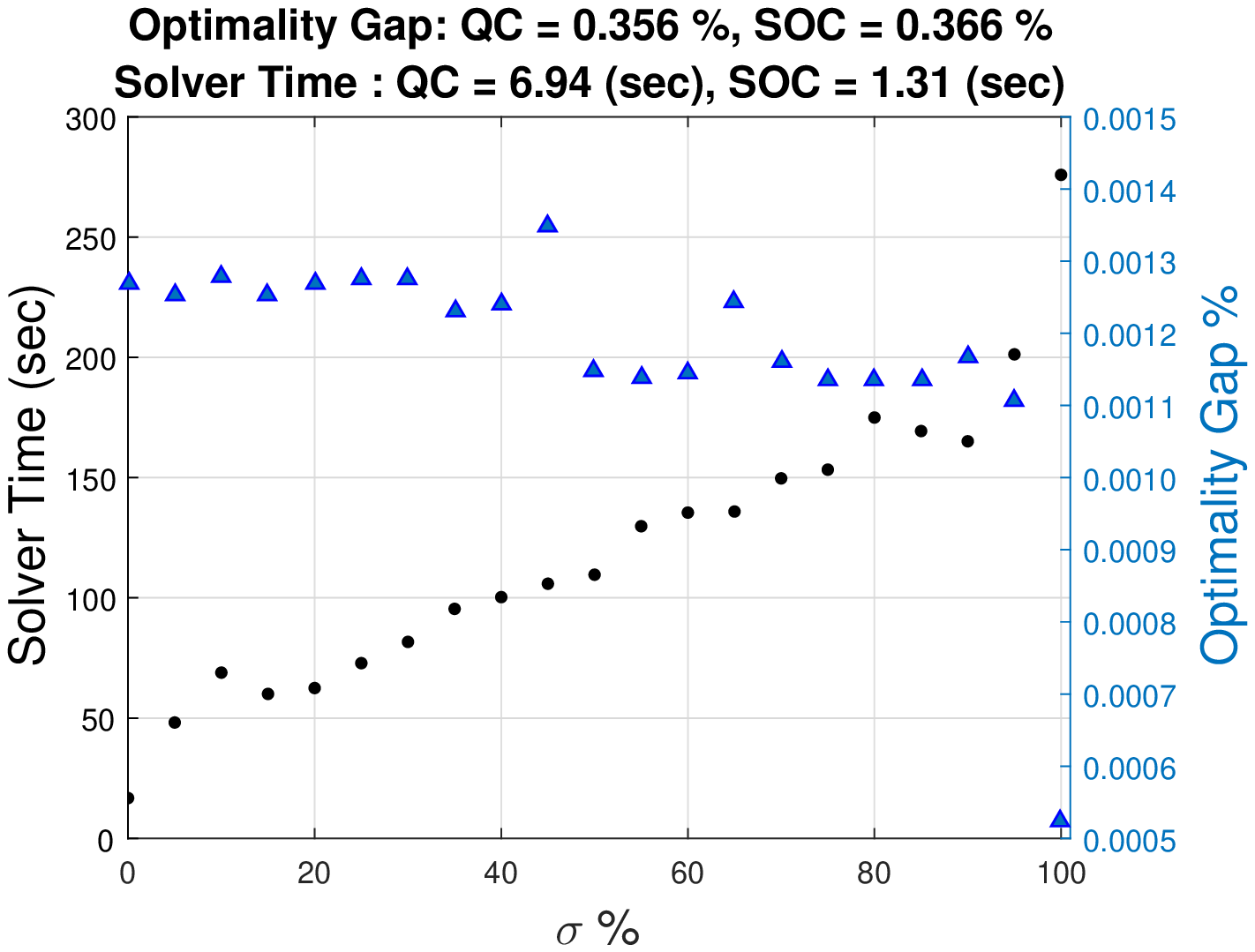}\label{case2746wop}}\hspace{2.5mm}
\subfloat[case3012wp\_mp\_\_api]{\includegraphics[totalheight=0.17\textheight]{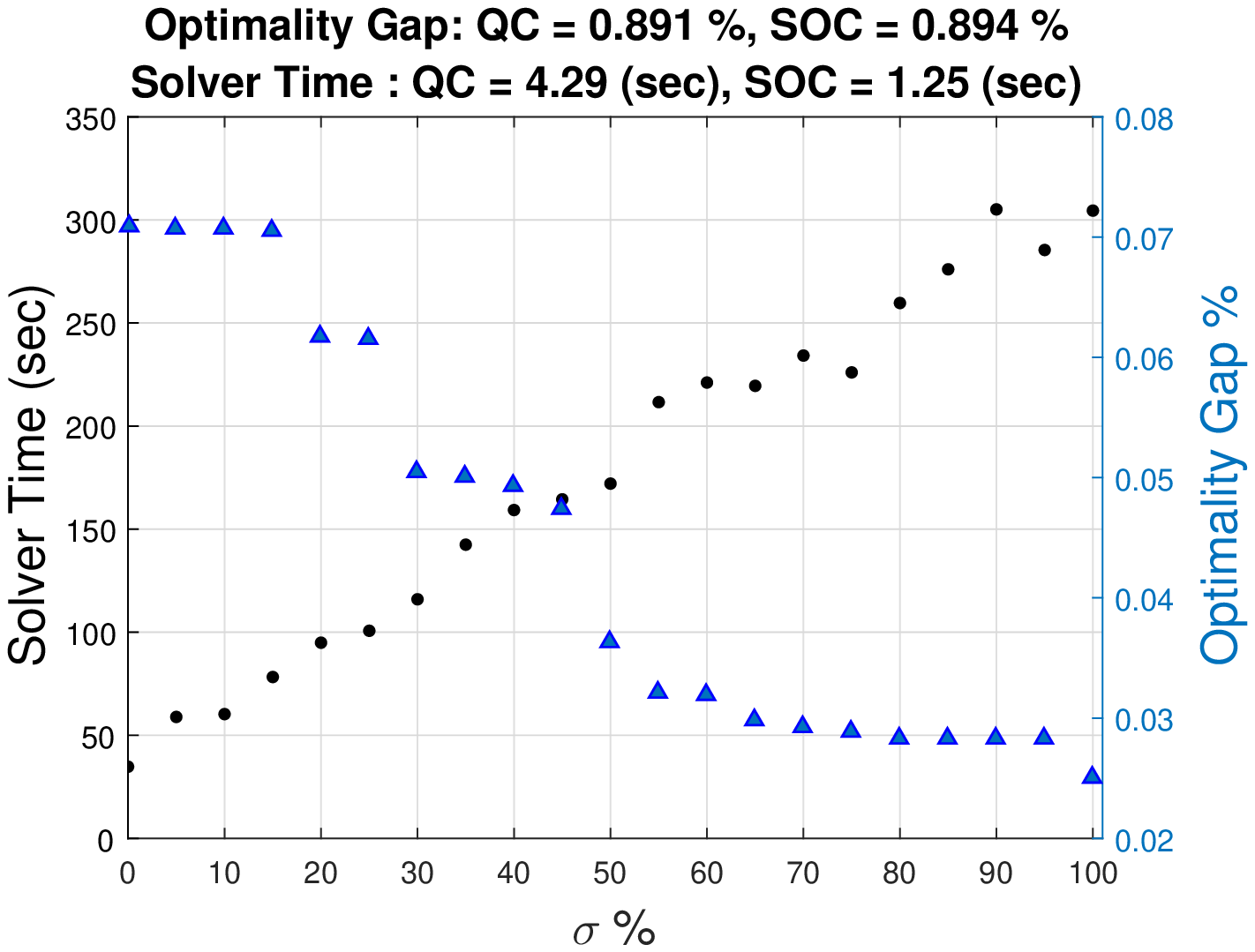}\label{3012wp_mp__api}}\hspace{2.5mm}
\subfloat[case2736sp\_mp\_\_api]{\includegraphics[totalheight=0.17\textheight]{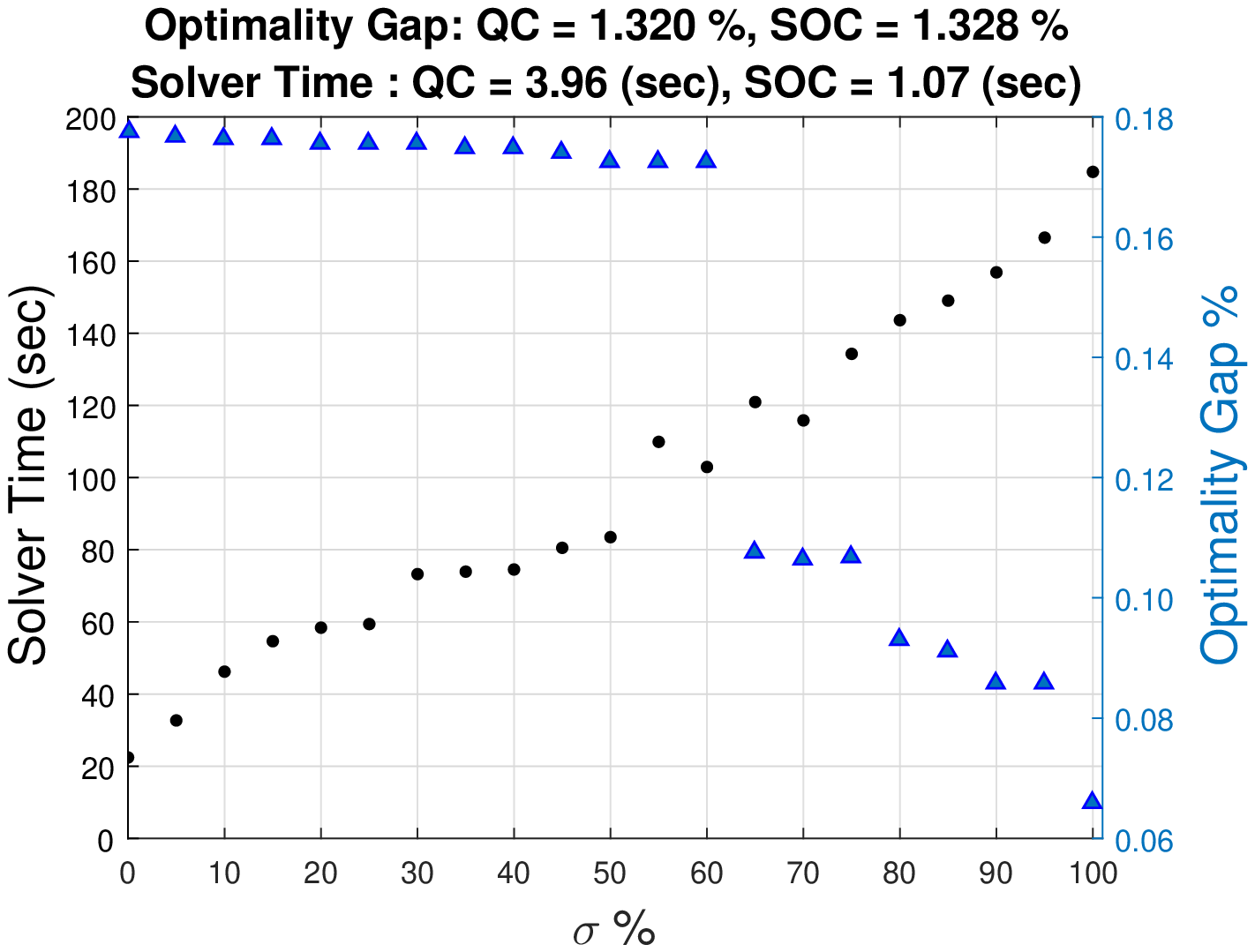}\label{2736sp_mp__api}}
\\
\subfloat[case2383wp\_mp]{\includegraphics[totalheight=0.17\textheight]{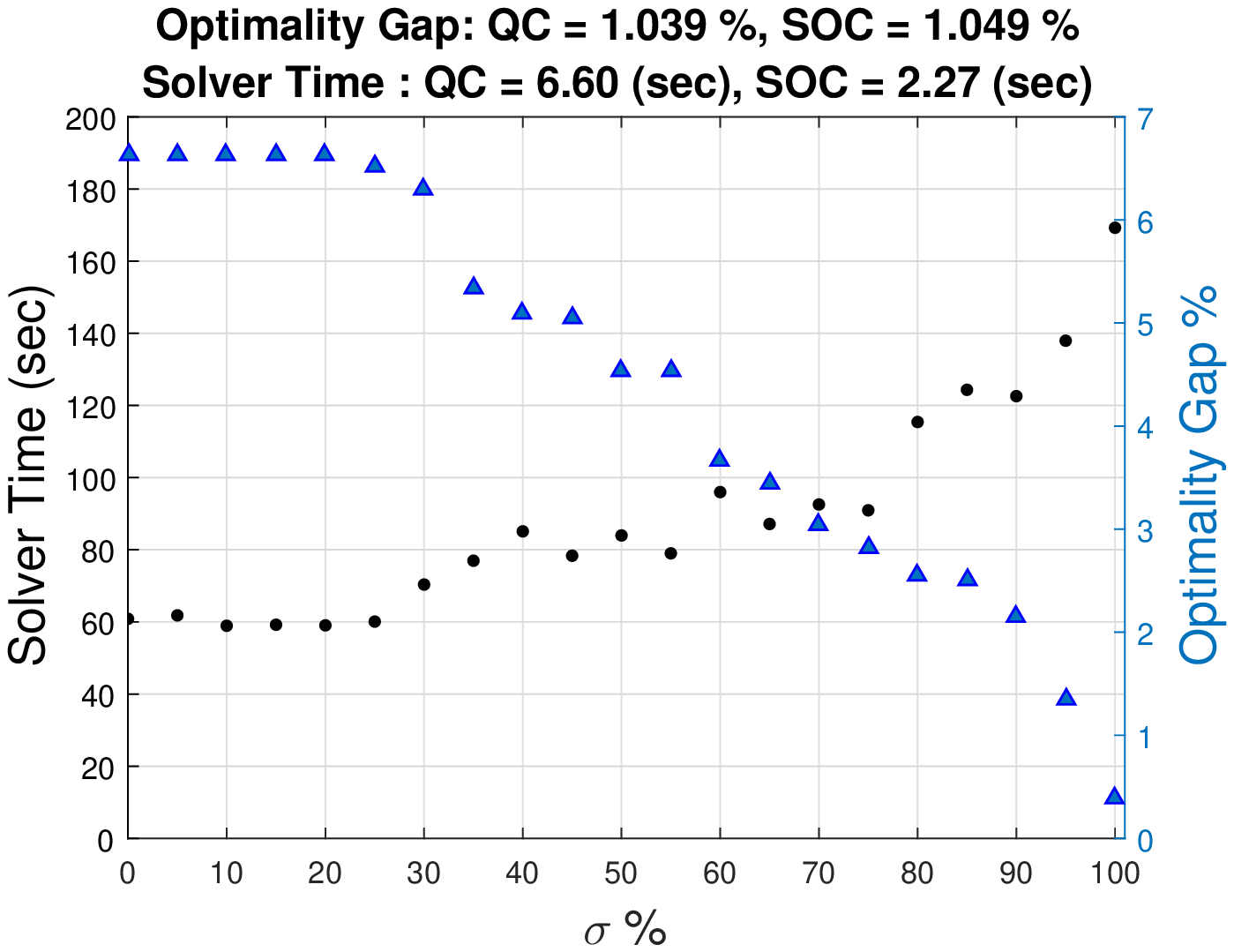}\label{2383wp}}\hspace*{2.5pt}
\subfloat[case3120sp\_mp\_\_api]{\includegraphics[totalheight=0.17\textheight]{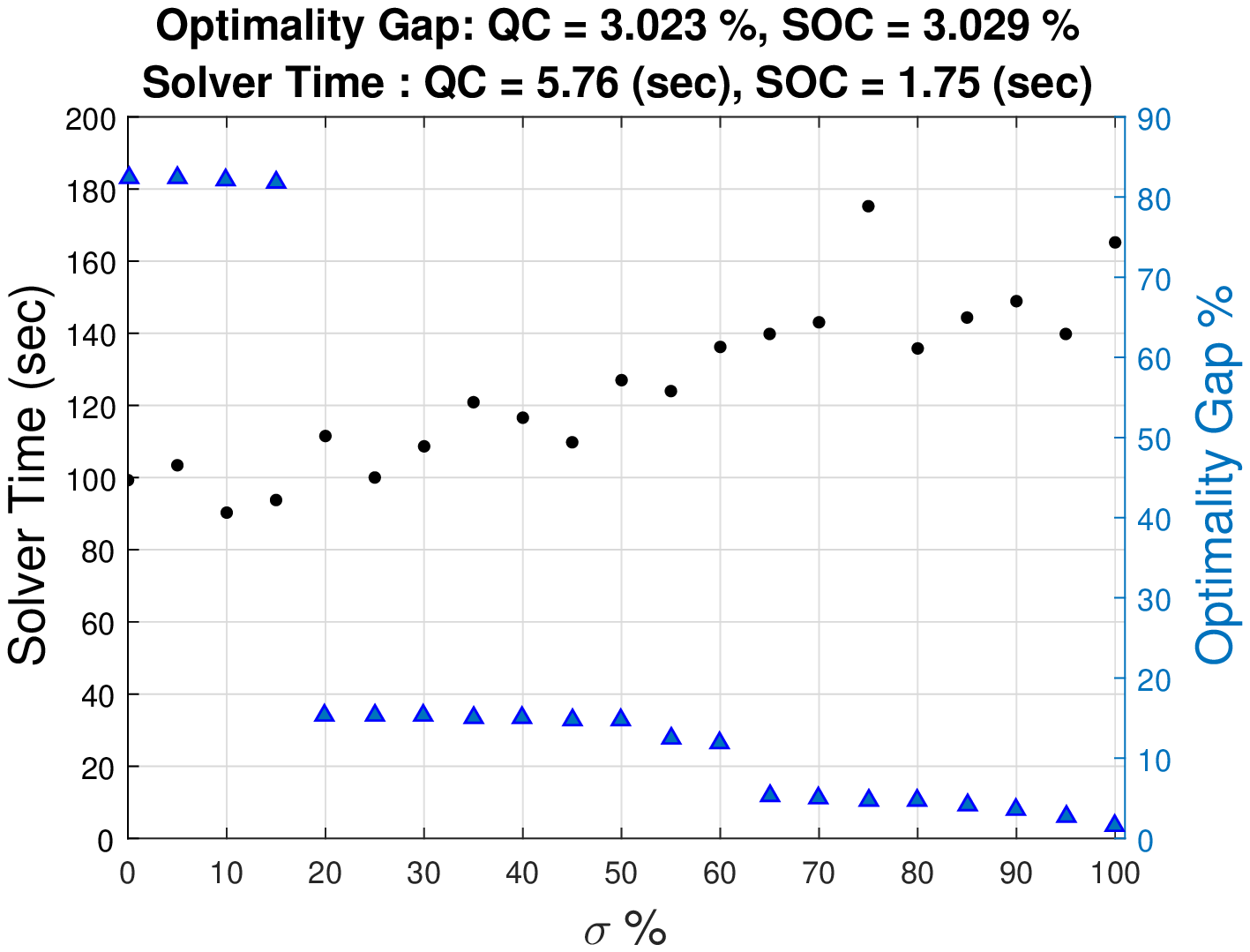}\label{3120sp_mp__api}}\hspace{2.5mm}
\subfloat[case3375wp\_mp]{\includegraphics[totalheight=0.17\textheight]{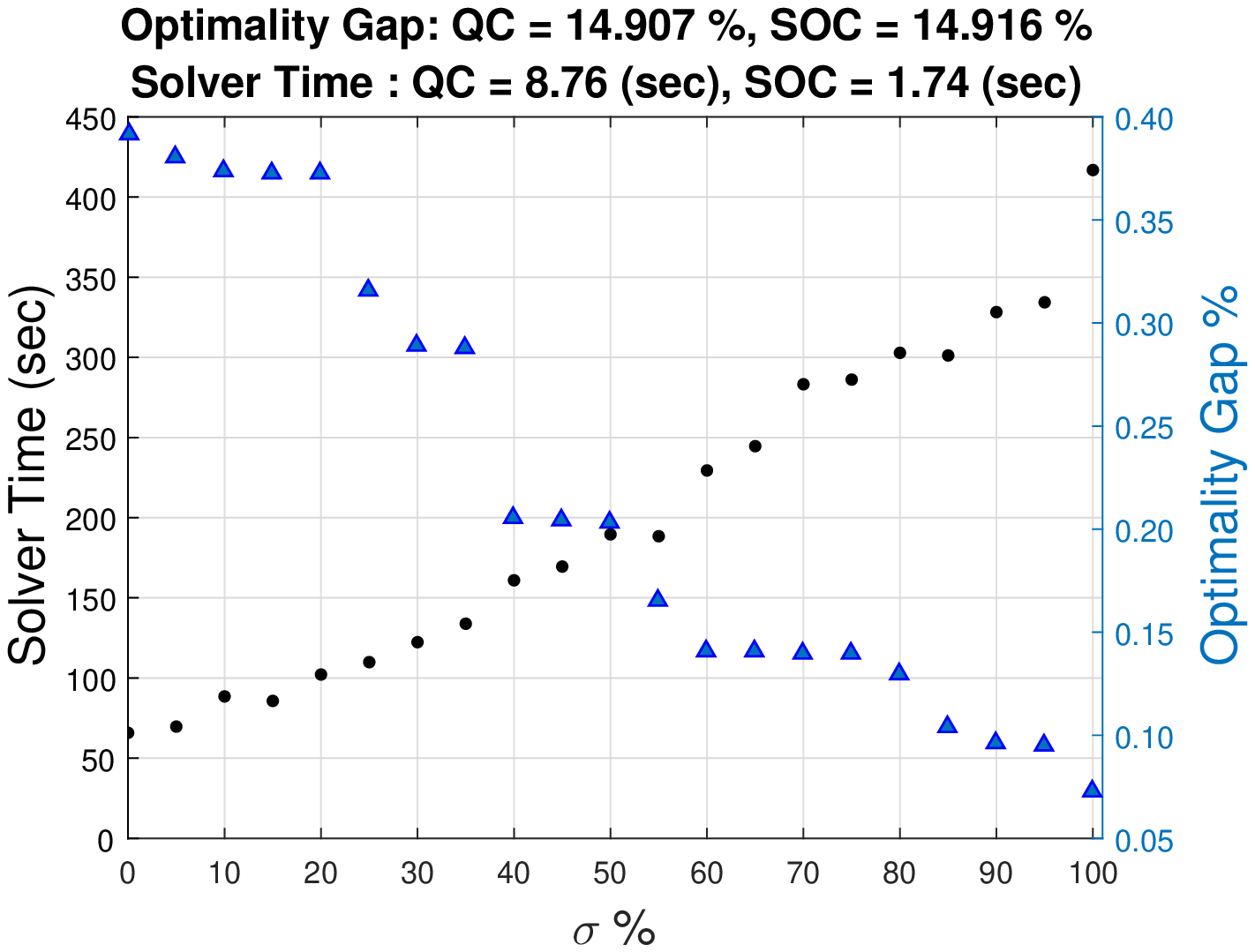}\label{3375wp}}
\caption{Solver Times (seconds) ($\bullet$) and Optimality Gaps (\%) ($\blacktriangle$) vs. $\sigma (\%)$.}
\label{Solver Times and Optimality Gaps}
\vspace*{-1em}
\end{figure*}

\begin{table*}[!t]
\centering
\caption{Optimality Gaps and Solver Times for $\sigma = 20\%$}
\label{Optimality Gaps_and_Solver_Times}
\begin{tabular}{|l||c||c||c|c|c|c||c|c|c|c|}
\hline
\multicolumn{1}{|c||}{\multirow{2}{*}{\begin{tabular}[c]{@{}c@{}}System \\ Model\end{tabular}}} & \multirow{2}{*}{\begin{tabular}[c]{@{}c@{}}MATPOWER \\ Obj. Value (\$/h)\end{tabular}} & \multicolumn{1}{c||}{\multirow{2}{*}{\begin{tabular}[c]{@{}c@{}}$\sigma$ (\%)\end{tabular}}}  &\multicolumn{4}{c||}{Optimality Gap (\%)} & \multicolumn{4}{c|}{Solver Time (sec)}\\ \cline{4-11}
                &    &                   &   Our Method  &  SDP      &  QC       &   SOCP    &   Our Method  &  SDP      &  QC       &   SOCP   \\ \hline
 \multicolumn{11}{|c|} {Typical Operating Condition (TOC)}  \\  \cline{1-11}
 case118\_ieee   & $\hphantom{0}3718.64$  & $20$ & $0.862$  & $0.064$  & $1.568$  & $2.069$  & $\hphantom{0}0.33$  &  $\hphantom{0}\hphantom{0}0.65$   &  $0.19$  & $0.11$ \\ \hline
 case300\_ieee   & $16891.28$  & $22$ & $0.148$  & $0.078$  & $1.175$  & $1.179$  & $\hphantom{0}0.86$  & $\hphantom{0}\hphantom{0}2.80$   &  $0.56$  & $0.18$ \\ \hline
 case1354\_pegase & $74069.35$  & $20$ & $0.107$  & $0.010$  & $0.077$  & $0.081$ & $\hphantom{0}4.46$  & $\hphantom{0}13.58$   & $2.31$  & $0.52$ \\ \hline
 case2383wp\_mp & $1868511.83\hphantom{0}\hphantom{0}$  & $22$ &  $6.643$  & $0.375$  & $1.039$  & $1.049$ & $58.48$ & $181.81$  & $6.60$  & $2.27$  \\ \hline
 case2736sp\_mp & $1307883.13\hphantom{0}\hphantom{0}$  & $20$ &  $0.000$  & $0.000$  & $0.288$  & $0.297$  & $50.00$  & $277.30$  & $6.29$  & $1.29$  \\ \hline
 case2737sop\_mp  & $777629.30\hphantom{0}$  & $20$ &  $0.001$  & $0.000$  & $0.246$  & $0.249$ & $53.10$  & $230.97$  & $3.49$  & $1.43$ \\ \hline
 case2746wop\_mp  & $1208279.81\hphantom{0}\hphantom{0}$  & $20$ &  $0.001$  & $0.001$  & $0.356$  & $0.366$ & $55.45$  & $318.05$  & $6.94$  & $1.31$ \\ \hline
 case2746wp\_mp   & $1631775.10\hphantom{0}\hphantom{0}$  & $20$ &  $0.000$  & $0.000$  & $0.318$  & $0.324$ & $51.55$  & $239.08$  & $6.21$  & $1.34$ \\ \hline
 case2869\_pegase & $133999.29\hphantom{0}$  & $20$ &  $0.085$  & $0.008$  & $0.088$  & $0.091$ & $13.02$  & $\hphantom{0}41.13$ & $4.52$  & $1.82$ \\ \hline
 case3012wp\_mp  & $2600842.77\hphantom{0}\hphantom{0}$  & $20$ &  $0.850$  & $0.143$  & $1.002$  & $1.021$ & $74.01$  & $362.15$  & $6.97$  & $1.27$ \\ \hline
 case3120sp\_mp & $2145739.43\hphantom{0}\hphantom{0}$  & $20$ &  $0.153$  & $0.088$  & $0.536$  & $0.545$ & $95.78$  & $435.99$  & $6.52$  & $1.34$  \\ \hline
 case3375wp\_mp & $7435697.54\hphantom{0}\hphantom{0}$  & $20$ &  $0.372$  & $0.073$  & $14.907\hphantom{0}$  & $14.916\hphantom{0}$  & $103.92\hphantom{0}$  & $449.06$  & $8.76$  & $1.74$  \\ \hline

 \multicolumn{11}{|c|} {Congested Operating Conditions (API)}  \\  \cline{1-11}
 case118\_ieee\_\_api  & $10325.27$  & $60$ & $37.888\hphantom{0}$  & $31.502\hphantom{0}$  & $43.650$\hphantom{0}  & $44.082$\hphantom{0} & $\hphantom{0}0.62$ & $\hphantom{0}\hphantom{0}0.84$  & $0.35$  & $0.14$  \\ \hline
 case300\_ieee\_\_api   & $22866.01$  &   $20$ &  $0.004$  & $0.003$  & $0.818$  & $0.837$  & $\hphantom{0}0.68$  & $\hphantom{0}\hphantom{0}1.97$  & $0.66$  & $0.23$  \\ \hline
 case2383wp\_mp\_\_api  & $23499.48$  &  $22$ &  $0.547$  & $0.101$  & $1.118$  & $1.120$ & $97.02$  & $194.25$  & $3.80$  & $1.15$ \\ \hline
 case2736sp\_mp\_\_api  & $25437.70$  &  $20$ &  $0.175$  & $0.069$  & $1.320$  & $1.328$ & $55.25$  & $191.73$  & $3.96$  & $1.07$ \\ \hline
 case2737sop\_mp\_\_api & $21192.40$  &  $20$ &  $0.318$  & $0.007$  & $1.049$  & $1.056$ & $41.77$  & $196.74$  & $3.94$  & $1.13$ \\ \hline
 case2746wop\_mp\_\_api & $22814.86$  &  $20$ &  $0.000$  & $0.000$  & $0.488$  & $0.491$ & $65.85$  & $282.61$ & $4.34$  & $1.40$ \\ \hline
 case2746wp\_mp\_\_api  & $27291.58$  &  $20$ &  $0.000$  & $0.000$  & $0.569$  & $0.491$ & $53.99$  & $199.81$  & $4.02$  & $1.15$ \\ \hline
 case3012wp\_mp\_\_api  & $27917.36$  &  $20$ &  $0.061$  & $0.025$  & $0.891$  & $0.894$ & $77.84$  & $266.51$  & $4.29$  & $1.25$  \\ \hline
 case3120sp\_mp\_\_api  & $22874.98$  &  $20$ &  $15.416$\hphantom{0}  & $0.542$  & $3.023$  & $3.029$ & $104.09$  & $458.98$  & $5.76$  & $1.75$  \\ \hline
 case3375wp\_mp\_\_api  & $48898.95$  &  $20$ &  $0.035$  & $0.014$  & $45.765\hphantom{0}$ & $45.766\hphantom{0}$  & $89.74$ & $284.44$ & $5.04$  & $1.63$  \\ \hline

 \multicolumn{11}{|c|} {Small Angle Difference Condition (SAD)}  \\  \cline{1-11}
 case118\_ieee\_\_sad & $\hphantom{0}4324.17$  &  $60$ &  $10.687\hphantom{0}$  & $7.652$  & $8.286$  & $15.783\hphantom{0}$ & $\hphantom{0}0.68$  & $\hphantom{0}\hphantom{0}0.93$  & $0.43$  & $0.17$ \\ \hline
 case2383wp\_mp\_\_sad& $1935308.26\hphantom{0}\hphantom{0}$  &  $31$ &  $22.467\hphantom{0}$  & $1.259$  & $2.967$  & $4.464$ & $82.19$  & $224.38$ & $6.66$  & $1.73$ \\ \hline
 case2736sp\_mp\_\_sad  & $1337042.83\hphantom{0}\hphantom{0}$  &  $20$ &  $38.915\hphantom{0}$  & $0.635$  & $2.009$  & $2.472$ & $66.92$  & $314.64$  & $4.39$  & $1.18$ \\ \hline
 case2746wop\_mp\_\_sad  & $1241955.37\hphantom{0}\hphantom{0}$  &  $20$ &  $20.126\hphantom{0}$  & $0.999$  & $2.482$  & $3.067$ & $78.31$  & $458.02$  & $5.25$  & $1.13$  \\ \hline
 case3012wp\_mp\_\_sad  & $2635451.43\hphantom{0}\hphantom{0}$  &  $20$ &  $15.537\hphantom{0}$  & $0.464$  & $1.916$ & $2.320$  & $103.00$ & $439.56$ & $5.19$  & $1.34$  \\ \hline
\end{tabular}
\vspace*{-1em}
\end{table*}

We evaluate bound tightness via the \emph{optimality gap}, which compares the objective value from the local solution to the bound from a dual feasible point obtained using the proposed method or the bound from a relaxation:
\begin{align*}
\small
&\text{\small Optimality Gap (\%)} = \frac{\text{\small Local Sol. Obj. -- Obj. Value Bound}}{\text{\small Local Sol. Obj.}} \times 100.
\end{align*}

\subsection{Test Case Results}
To illustrate the trade-off in computational speed and bound tightness, Fig.~\ref{Solver Times and Optimality Gaps} shows how the solver times and optimality gaps vary with $\sigma$ for a representative selection of test cases. For some cases, such as case3012wp\_mp\_\_api, case2736sp\_mp\_\_api, and case3375wp\_mp, the solver times decrease by approximately a factor of nine to ten as $\sigma$ varies from $100\%$ (equivalent to the sparse version of the SDP relaxation~\eqref{eq:Dualsdpopf}) to $0\%$ (only submatrices that are not positive semidefinite are identified as problematic). Varying $\sigma$ has a much smaller impact on the solver times for other cases, such as case2746wop\_mp, case2383wp\_mp, and case3120sp\_mp\_\_api.

The improvement in speed with smaller $\sigma$ has a trade-off with respect to the tightness of the objective value bound. As with the solver times, the impact of $\sigma$ on the optimality gap is case dependent. For case2746wop\_mp, case3012wp\_mp\_\_api, and case2736sp\_mp\_\_api, the optimality gaps decrease by approximately a factor of two to three as $\sigma$ varies from $0\%$ to $100\%$, while for case3375wp\_mp, case2383wp\_mp, and case3120sp\_mp\_\_api, the optimality gaps decrease by factors of 8, 20 and 200, respectively, as $\sigma$ varies from $0\%$ to $100\%$. These variations are explained by the quality of the local solution. When the local solution is close to the SDP solution, increasing $\sigma$ only has a small impact on the bound.

In contrast to the solver times, which vary gradually with $\sigma$, the optimality gaps exhibit more significant variations. For example, the optimality gap for case2736sp\_mp\_\_api  does not change significantly from $\sigma=0\%$ until a steep decrease at $\sigma=60\%$. This is explained by a group of buses whose dual variable values from the local solution are far from their values in the SDP solution. Choosing $\sigma < 60$\% fixes these dual variables to the values from the local solution, restricting our method's ability to achieve a tight bound. Jumps in the optimality gaps for other test cases have similar explanations.

Table~\ref{Optimality Gaps_and_Solver_Times} shows results for representative test cases, specifically comparing the proposed method with the SDP, QC, and SOCP relaxations. Three categories of NESTA test cases are considered: typical operating conditions (TOC), congested operation (API), and small angle difference (SAD). For these results, we choose $\sigma$ to be either $20\%$ or the lower bound on $\sigma$ needed to include all non-positive-semidefinite submatrices as discussed in Step 3 of our proposed algorithm in Section~\ref{Algo_descrip}, whichever is greater for each test case. Numerical results, including studies not detailed here, suggest that this value for $\sigma$ yields good performance for a variety of test cases. Further tuning of $\sigma$ can improve the results achieved for specific problems.

Table~\ref{Average_Optimality_Gaps_and_Solver_Times} presents the average optimality gaps and solver times for each category of test cases. The averages in this table are computed based on the values in Table~\ref{Optimality Gaps_and_Solver_Times} using the values of $\sigma$ reported therein. Our method performs well for most of the TOC and API test cases, as indicated by smaller optimality gaps than the QC and SOCP relaxations with a significant (approximately factor of five) improvement in computation speed relative to the SDP relaxation. The average results show that the proposed method yields better optimality gaps than QC and SOCP relaxations, while being significantly faster than the SDP relaxation, suggesting the method's utility as a middle ground between the the SDP relaxation and the SOCP and QC relaxations.
%
%
\begin{table}[ht]
\vspace*{-0.75em}
\centering
\caption{Average Optimality Gaps and Solver Times}
\label{Average_Optimality_Gaps_and_Solver_Times}
\begin{tabular}{|l||c|c|c||c|c|c|}
\hline
\multicolumn{1}{|c||}{\multirow{2}{*}{Approach}} & \multicolumn{3}{c||}{\begin{tabular}[c]{@{}c@{}}Average\\ Optimality Gaps (\%)\end{tabular}} & \multicolumn{3}{c|}{\begin{tabular}[c]{@{}c@{}}Average\\ Solver Times (sec)\end{tabular}} \\ \cline{2-7} 
      & TOC  & API    & SAD  & TOC   & API   & SAD     \\ \hline
Our Method    & \hphantom{0}0.77\hphantom{0}                        & 1.84   &  21.55\hphantom{0}                                    & \hphantom{0}46.9   & \hphantom{0}65.6   &  66.2   \\ \hline
SDP  & 0.07   & 0.08       &    2.20    & 212.7                       & 230.8     &   287.5\hphantom{0}    \\ \hline
QC      & 1.80        & 6.12        &       3.53                   & \hphantom{0}\hphantom{0}5.0  & \hphantom{0}\hphantom{0}4.0    &      \hphantom{0}4.4  \\ \hline
SOCP & 1.85  & 6.11     &  5.62  & \hphantom{0}\hphantom{0}1.2 & \hphantom{0}\hphantom{0}1.2     &   \hphantom{0}1.1  \\ \hline
\end{tabular}
\vspace*{-0.25em}
\end{table}

However, despite promising performance on many problems, there are test cases which challenge our method. For instance, our method yields large gaps relative to the SOCP, QC, and SDP relaxations for the SAD test cases, case2383wp\_mp, and case3120sp\_mp\_\_api. For these cases, the local solution is not close to the SDP relaxation's solution such that large values of $\sigma$ (with correspondingly slow solution times) are necessary to achieve small optimality gaps. Explaining our method's poor performance for these test cases, the dual correspondence ratios~\eqref{eq:norm_Dist} presented in Table~\ref{Dual_Correspondence_Ratio_and_Percentage_of_Positive_Semidefinite_Submatrices} have larger values for the test cases where our method performs poorly relative to the test cases with tighter optimality gaps.
%
%
\begin{table}[t]
\centering
\caption{Dual Correspondence Ratios and Percentages of Positive Semidefinite Submatrices for Selected Test Cases}
\label{Dual_Correspondence_Ratio_and_Percentage_of_Positive_Semidefinite_Submatrices}
\begin{tabular}{|l|c|c|c|}
\hline
\multicolumn{1}{|c|}{\begin{tabular}[c]{@{}c@{}}System\\ Model\end{tabular}} & \begin{tabular}[c]{@{}c@{}}Dual \\ Correspondence\\ Ratio (\%)\end{tabular} & \begin{tabular}[c]{@{}c@{}}Percentage \\ of PSD\\ Submatrices \end{tabular} & \begin{tabular}[c]{@{}c@{}}Optimality \\ Gaps (\%)\end{tabular}\\ \hline
case2736sp\_mp                      & \hphantom{0}\hphantom{0}0.16        & 99.3  &  0.000\\ \hline
case2746wop\_mp                     & \hphantom{0}\hphantom{0}0.45        & 98.4  & 0.001 \\ \hline
case3375wp\_mp\_\_api               & \hphantom{0}\hphantom{0}0.93        & 97.5  &  0.035\\ \hline
case3012wp\_mp                      & \hphantom{0}\hphantom{0}1.32        & 97.3  &  0.850\\ \hline
case2869\_pegase                    & \hphantom{0}\hphantom{0}4.32        & 96.7  &  0.085\\ \hline
case3375wp\_mp                      & \hphantom{0}\hphantom{0}7.00        & 95.7  &  0.372\\ \hline
case2383wp\_mp                      & \hphantom{0}28.91                   & 78.0   &  6.643\\ \hline
case3120sp\_mp\_\_api               & 462.70                              & 92.0    & 15.416\hphantom{0}\\ \hline
case2383wp\_mp\_\_sad               & 105.10    
        & 58.4    &  22.467\hphantom{0}   \\ \hline
case2736sp\_mp\_\_sad               & \hphantom{0}61.30                   &  89.8   & 38.915\hphantom{0}\\ \hline
case3012wp\_mp\_\_sad               & 104.31    
        &  86.2   &  15.537\hphantom{0} \\ \hline
\end{tabular}
\vspace*{-1em}
\end{table}

Finally, we demonstrate the a-priori heuristic proposed in Section~\ref{III:discustion} via selected test cases in Table~\ref{Dual_Correspondence_Ratio_and_Percentage_of_Positive_Semidefinite_Submatrices}. Our method performs well for the test cases where the percentages of positive semidefinite submatrices are greater than $95$\%. For these cases, our method yields optimality gaps between those from the SDP relaxation and the SOCP and QC relaxations with significant speed improvements compared to the SDP relaxation. Conversely, our method often performs poorly for test cases that do not satisfy this heuristic, such as case2383wp\_mp, case3120sp\_mp\_\_api, and the SAD cases. This shows the appropriateness of our heuristic recommendation for applying our method to problems where the percentage of submatrices $\hat{\mathbf{A}}_i$ that are positive semidefinite is greater than $95$\%.

\section{Conclusions}
\label{l:conclusion}

This paper has proposed an SDP-based method that leverages knowledge of a local solution to quickly compute objective value bounds for optimal power flow problems. Numerical experiments indicate that the proposed method is approximately four to five times faster than the SDP relaxation for a variety of test cases without excessive degradation in the quality of the objective bound. These experiments also show that the proposed method is slower than the SOCP and QC relaxations but can obtain tighter bounds. Thus, the proposed method provides a middle ground between the SDP relaxation and the SOCP and QC relaxations with a trade-off between computational speed and tightness.



\bibliographystyle{IEEEtran}
\bibliography{IEEEabrv,Quick_Calc_Lower_Bound}
\end{document}